\documentclass[10pt]{amsart}
\usepackage{amssymb, amsthm, amsmath}
\usepackage{xy}
\usepackage{epsfig}
\usepackage{psfrag}
\usepackage{xypic}

\newcommand{\tto}{\twoheadrightarrow}
\newcommand{\X}{{\mathfrak X}}
\newcommand{\Y}{{\mathcal Y}}
\newcommand{\Z}{{\mathcal Z}}
\newcommand{\C}{{\mathbb C}}
\newcommand{\G}{{\mathcal G}}

\newcommand{\Groth}{{\mathfrak{G}}}
\newcommand{\bull}{{\sssize \bullet}}
\newcommand{\rank}{\operatorname{rank}}
\newcommand{\Hom}{\operatorname{Hom}}
\newcommand{\Fl}{\operatorname{F\ell}}
\newcommand{\Gr}{\operatorname{Gr}}
\renewcommand{\O}{{\mathcal O}}
\newcommand{\wt}{\widetilde}
\newcommand{\gtr}{\vartriangleright}

\newtheorem{thm}{Theorem} 
\newtheorem*{thm1a}{Theorem 1$'$}
\newtheorem*{thmII}{Theorem 2}
\newtheorem{lemma}{Lemma} 
\newtheorem{prop}{Proposition}
\newtheorem{cor}{Corollary}

\theoremstyle{definition}
\newtheorem{example}{Example}

\newcommand{\refprop}[1]{Proposition~\ref{#1}}

\begin{document}

\title[Grothendieck polynomials and quiver formulas]
{Grothendieck polynomials and quiver formulas}
\author{Anders S.~Buch, Andrew Kresch, 
Harry Tamvakis, and Alexander Yong}
\date{June 20, 2003}
\thanks{The third author was supported in part by 
  NSF Grant DMS-0296023.}
\address{Matematisk Institut, Aarhus Universitet, Ny Munkegade, 8000
  {\AA}rhus C, Denmark}
\email{abuch@imf.au.dk}
\address{Department of Mathematics, University of Pennsylvania,
209 South 33rd Street,
Philadelphia, PA 19104-6395, USA}
\email{kresch@math.upenn.edu}
\address{Department of Mathematics, Brandeis University - MS 050,
P. O. Box 9110, Waltham, MA
02454-9110, USA}
\email{harryt@brandeis.edu}
\address{Department of Mathematics,  
University of Michigan, 525 East University Ave.,
Ann Arbor, MI 48109-1109, USA}
\email{ayong@umich.edu}

\begin{abstract}
  Fulton's universal Schubert polynomials give cohomology formulas for
  a class of degeneracy loci, which generalize Schubert varieties.
  The $K$-theoretic quiver formula of Buch expresses the structure
  sheaves of these loci as integral linear combinations of products of
  stable Grothendieck polynomials.  We prove an explicit combinatorial
  formula for the coefficients, which shows that they have alternating
  signs.  Our result is applied to obtain new expansions for the
  Grothendieck polynomials of Lascoux and Sch\"utzenberger.
\end{abstract}

\maketitle

%%%%%%%%%%%%%%%%%%%%%%%%%%%%%%%%%%%%%%%%%%%%%%%%%%%%%%%%%%%%%%%%%%%%%%
%%%   1 .  INTRODUCTION
%%%%%%%%%%%%%%%%%%%%%%%%%%%%%%%%%%%%%%%%%%%%%%%%%%%%%%%%%%%%%%%%%%%%%%

\section{Introduction and main results}

Let $\X$ be a smooth complex algebraic variety and let
\begin{equation} \label{eqn:univseq}
E_1 \to \cdots \to E_{n-1}\to E_n \to F_n \to F_{n-1} \to \cdots \to F_1
\end{equation}
be a sequence of vector bundles and morphisms over $\X$, such that
$\rank(F_i) = \rank(E_i) = i$ for $1 \leq i \leq n$. For any
permutation $w \in S_{n+1}$, there is a degeneracy locus
\begin{equation} \label{eqn:univlocus}
  \Omega_w(E_\bull \to F_\bull) = 
  \{ x \in \X \mid \rank(E_q(x) \to F_p(x)) \leq r_w(p,q), ~\forall 
  \ 1\leq p,q\leq n \} \,,
\end{equation}
where $r_{w}(p,q)$ is the number of $i\leq p$ such that $w(i)\leq q$.  We will
assume that the bundle maps are sufficiently general so that this degeneracy
locus has the expected codimension, equal to the length $\ell(w)$. In this
situation, Fulton \cite{fulton:usp} gave a formula for the cohomology class of
$\Omega_{w} = \Omega_{w}(E_\bull\to F_\bull)$ in $H^{*}(\X,{\mathbb Z})$ as a
{\em universal Schubert polynomial} in the Chern classes of the vector bundles
involved.

While the cohomology class of $\Omega_{w}$ gives useful global
information, there is even more data hidden in its structure sheaf
${\mathcal O}_{\Omega_{w}}$. The main result of this paper gives an
explicit combinatorial formula for the class $[\O_{\Omega_w}]$ of this
structure sheaf in the Grothendieck ring $K(\X)$ of algebraic vector
bundles on $\X$.  To state it, we need the {\em degenerate Hecke
  algebra}, which is the associative ${\mathbb Z}$-algebra
generated by $s_i$ for $i = 1, 2, \dots$ with relations
\begin{eqnarray*}
s_{i}^2 & = & s_{i} \\
s_i s_j & = & s_j s_i \text{\ \ \ \ \ for $|i-j| > 1$} \\
s_{i} s_{i+1} s_{i} & = & s_{i+1} s_{i} s_{i+1} \,.
\end{eqnarray*} 
We also require the stable Grothendieck polynomials $G_{u}(E - E')$, where $u$
is a permutation and $E$, $E'$ are vector bundles over $\X$ (see Section
\ref{prelims} for the definition).

\begin{thm} \label{thm:main1}
  For $w \in S_{n+1}$ we have
\[ [\O_{\Omega_w}] = 
   \sum (-1)^{\ell(u_1 \cdots u_{2n-1} w)} \, G_{u_1}(E_2-E_1)
   \cdots G_{u_n}(F_n- E_n) \cdots G_{u_{2n-1}}(F_1- F_2)
\]
in $K(\X)$, where the sum is over all factorizations $w = u_1 \cdots u_{2n-1}$
in the degenerate Hecke algebra such that $u_i \in S_{\min(i,2n-i)+1}$ for
each $i$.
\end{thm}

The above formula corresponds to computing the alternating sum of a
locally free resolution of ${\mathcal O}_{\Omega_{w}}$ in $K(\X)$, and
thus includes a formula for the cohomology class of $\Omega_w$ as its
leading term. Theorem \ref{thm:main1} is therefore a generalization
of \cite[Thm.~3]{schquiv}.

The locus $\Omega_w(E_\bull \to F_\bull)$ is a special case of a {\em
  quiver variety}.  In \cite{buch:grothendieck} a formula for the
  class of the structure sheaf of a general quiver variety is proved,
  which expresses this class as a linear combination of products of
  stable Grothendieck polynomials for Grassmannian permutations.
  Furthermore, it is conjectured that the {\em quiver coefficients}
  occurring in this formula have signs which alternate with the
  codimension.

The quiver formula specializes to {\em
universal Grothendieck polynomials\/} $\G_w(F_\bull; E_\bull)$ in the
exterior powers of the bundles (and the inverse of the top powers),
which are $K$-theoretic analogues of universal Schubert polynomials.  
Given any partition $\alpha$, we let $G_\alpha = G_{w_\alpha}$ denote
the stable Grothendieck polynomial for the Grassmannian permutation 
$w_{\alpha}$ corresponding to $\alpha$.  Then the quiver formula 
has the form
\[ [\O_{\Omega_w}]=\G_w(F_\bull; E_\bull) = \sum_\lambda c_{w,\lambda}^{(n)} \,
   G_{\lambda^1}(E_2-E_1) \cdots G_{\lambda_n}(F_n-E_n) \cdots
   G_{\lambda^{2n-1}}(F_1-F_2)
\]
where the sum is over finitely many sequences of partitions $\lambda =
(\lambda^1, \dots, \lambda^{2n-1})$ and the $c_{w,\lambda}^{(n)}$ are
quiver coefficients.

Theorem \ref{thm:main1} combined with a result of Lascoux
\cite{lascoux:grothendieck} proves that these coefficients do in fact have
alternating signs. Define integers $a_{w,\beta}$ such that $G_w = \sum
a_{w,\beta} G_\beta$, the sum over all partitions $\beta$.  Lascoux has shown
that $a_{w,\beta}$ is equal to $(-1)^{|\beta|-\ell(w)}$ times the number of
paths from $w$ to $w_\beta$ in a graph of permutations.  Given this result,
Theorem \ref{thm:main1} is equivalent to the following explicit combinatorial
formula for quiver coefficients:
\[
  c_{w,\lambda}^{(n)} = (-1)^{|\lambda| - \ell(w)}
  \sum_{u_1\cdots u_{2n-1} = w} 
  |\, a_{u_1,\lambda^1} \, a_{u_2,\lambda^2} \cdots
  a_{u_{2n-1},\lambda^{2n-1}} | \,.
\]

Our proof of Theorem \ref{thm:main1} is based on a special case of
this formula, proved in \cite{buch:grothendieck}, together with the
following Cauchy identity, which provides a $K$-theoretic
generalization of \cite[Cor.\ 2]{kirillov:cauchy} (see also
\cite[Thm.~3.7]{fulton:usp}).

\begin{thm}[Cauchy formula]
\label{thm:main2}
Let $E_{i}, F_{i}$, and $H_{i}$ for $i=1,\ldots, n$ be three collections
of vector bundles on $\X$. Then for any $w\in S_{n+1}$, we have
\[ \G_{w}(F_\bull; E_\bull) = \sum_{u\cdot v=w} (-1)^{\ell(uvw)} \,
   \G_{u}(H_\bull; E_\bull) \, \G_{v}(F_\bull; H_\bull)
\]
where the sum is over all permutations $u,v$ such that the product of
$u$ and $v$ is equal to $w$ in the degenerate Hecke algebra.
\end{thm} 

As a further consequence of our results, we obtain new formulas for
the double Grothendieck polynomials of Lascoux and Sch\"utzenberger
\cite{lascoux.schutzenberger:grothendieck},
which express these polynomials as linear combinations of stable
Grothendieck polynomials in disjoint intervals of variables.  The
coefficients in these expansions are all quiver coefficients; in
particular, this is true for the monomial coefficients of Grothendieck
polynomials.

This paper is organized as follows.  We review the facts about
Grothendieck polynomials that we require in Section \ref{prelims}. The
quiver varieties and universal Grothendieck polynomials are introduced
in Section \ref{UGP}. We prove the Cauchy formula in Section \ref{CF},
while our main theorem is proved in the following section. Finally in
Section \ref{SGP} we apply our results to obtain splitting formulas
for double Grothendieck polynomials.

The third author wishes to thank Marc Levine and the Universit\"at 
Essen for their hospitality and the Wolfgang Paul program of the 
Humbolt Foundation for support during the latter stages of work on 
this article.

%%%%%%%%%%%%%%%%%%%%%%%%%%%%%%%%%%%%%%%%%%%%%%%%%%%%%%%%%%%%%%%%%%%%%%
%%%   2 .  GROTHENDIECK POLYNOMIALS
%%%%%%%%%%%%%%%%%%%%%%%%%%%%%%%%%%%%%%%%%%%%%%%%%%%%%%%%%%%%%%%%%%%%%%

\section{Grothendieck polynomials}
\label{prelims}

We begin by recalling the definition of Lascoux and Sch{\"u}tzenberger's 
double Grothendieck polynomials \cite{lascoux.schutzenberger:grothendieck}.
Let $X=(x_1,x_2,\ldots)$ and $Y=(y_1,y_2,\ldots)$ be two sequences of
commuting independent variables and $w\in S_{n}$. If
$w = w_0$ is the longest permutation in $S_n$, then we set 
\[ \Groth_{w_0}(X;Y) = \prod_{i+j\leq n} (x_i+y_j-x_i y_j) \,. \]
If $w\neq w_0$, we can find a simple transposition $s_i=(i, i+1)\in
S_{n}$ such that $\ell(ws_i)=\ell(w)+1$.  We then define
\[ \Groth_w = \pi_i(\Groth_{ws_i}) \]
where $\pi_i$ is the isobaric divided difference operator given by
\[ \pi_i(f) = \frac{(1-x_{i+1})f(x_1,\dots,x_i,x_{i+1},\ldots,x_n) - 
   (1-x_i)f(x_1,\ldots,x_{i+1},x_i,\ldots,x_n)}{x_i-x_{i+1}} \,.
\]

Given permutations $u_1,\dots,u_m$, and $w$, we will write $u_1 \cdots
u_m = w$ if the product of the $u_i$ is equal to $w$ in the degenerate
Hecke algebra.  With this notation, the Grothendieck polynomials
satisfy the following Cauchy identity, which is due to Fomin and
Kirillov (see \cite[Thm.~8.1]{fomin.kirillov:yang-baxter} and
\cite{fomin.kirillov:grothendieck}):
\begin{equation}
\label{eqn:cauchy_old}
  \Groth_{w}(X;Y) = \sum_{u\cdot v=w}
  (-1)^{\ell(uvw)}\, \Groth_{u}(0;Y)\, \Groth_{v}(X;0) \,.
\end{equation}

Next we recall the definition of stable Grothendieck polynomials.
Given a permutation $w\in S_{n}$, and a nonnegative integer $m$, let
$1^{m}\times w\in S_{m+n}$ denote the shifted permutation which is the
identity on $\{1,2,\ldots,m\}$ and which maps $j$ to $w(j-m)+m$ for
$j>m$. It is known \cite{fomin.kirillov:grothendieck, 
fomin.kirillov:yang-baxter} that when $m$ grows to infinity, the
coefficient of each fixed monomial in $\Groth_{1^{m}\times w}$
eventually becomes stable. The double stable Grothendieck polynomial
$G_w \in {\mathbb Z}[[X;Y]]$ is the resulting power series:
\[ G_w = G_w(X;Y) = \lim_{m\to\infty}\Groth_{1^m\times w}(X;Y) \,. \]
The power series $G_w(X;Y)$ is symmetric
in the $X$ and $Y$ variables separately, and 
\[ G_w(1-e^{-X};1-e^{Y}) = 
   G_w(1-e^{-x_1},1-e^{-x_2},\ldots; 1-e^{y_1},1-e^{y_2},\ldots)
\]
is super-symmetric, that is, if one sets $x_1 = y_1$ in this expression,
then the result is independent of $x_1$ and $y_1$.

In particular, we will need stable Grothendieck polynomials for Grassmannian
permutations. If $\alpha=(\alpha_1\geq \alpha_2 \geq \alpha_3\geq \cdots)$ is
a partition and $p\geq \ell(\alpha)$, i.e., $\alpha_{p+1}=0$, the Grassmannian
permutation for $\alpha$ with descent in position $p$ is the permutation
$w_{\alpha}$ such that $w_{\alpha}(i)=i+\alpha_{p+1-i}$ for $1\leq i\leq p$
and $w_{\alpha}(i)< w_{\alpha}(i+1)$ for $i\neq p$. Now let
$G_{\alpha}=G_{w_{\alpha}}$; this is independent of the choice of $p$.
According to \cite[Thm.~6.6]{buch:littlewood-richardson} there are integers
$d_{\beta\gamma}^{\alpha}$ (with alternating signs) such that
\begin{equation}
\label{refequ}
G_{\alpha}(X;Y) = 
\sum d_{\beta\gamma}^{\alpha}\, G_{\beta}(X;0)\, G_{\gamma}(0;Y) \,.
\end{equation}

%Notice that if $q>p$ then the
%Grassmannian permutation for $\lambda$ with descent at position $q$ is
%equal to $1^{q-p}\times w_{\lambda}$. Thus $G_{\lambda}$ is
%independent of the choice of $p$.

Let $\Gamma\subseteq {\mathbb Z}[[X;Y]]$ be the linear span of all
stable Grothendieck polynomials.  It is shown in
\cite{buch:littlewood-richardson} that $\Gamma$ is closed under
multiplication and that the elements $G_{\alpha}$ form a ${\mathbb
  Z}$-basis of $\Gamma$.  In fact $\Gamma$ is a commutative and
cocommutative bialgebra with the coproduct $\Delta : \Gamma \to
\Gamma\otimes\Gamma$ given by $\Delta G_{\alpha} = \sum_{\beta,
  \gamma} d^{\alpha}_{\beta\gamma}\, G_{\beta} \otimes G_{\gamma}$.

We will describe a formula of Lascoux for the expansion of a stable
Grothendieck polynomial $G_w$ as a linear combination of these
elements.  Let $r$ be the last descent position of $w$, {i.e.\ }$r$ is
maximal such that $w(r) > w(r+1)$.  Set $w' = w \tau_{rk}$ where $k >
r$ is maximal such that $w(r) > w(k)$.  We also set $I(w) = \{ i < r
\mid \ell(w' \tau_{ir}) = \ell(w) \}$.

Define a relation $\gtr$ on the set of all permutations as follows.
If $I(w) = \emptyset$ we write $w \gtr v$ if and only if $v = 1 \times
w$.  Otherwise we write $w \gtr v$ if and only if there exist elements
$i_1 < \dots < i_p$ of $I(w)$, $p \geq 1$, such that $v = w' \tau_{i_1
  r} \dots \tau_{i_p r}$.  The following is an immediate consequence
of \cite[Thm.~4]{lascoux:grothendieck}.

\begin{thm}[Lascoux]
\label{thm:lascoux}
  For any permutation $w$ we have
\[ G_w = \sum_\beta a_{w,\beta} \, G_\beta \]
where the sum is over all partitions $\beta$, and $a_{w, \beta}$ is
equal to $(-1)^{|\beta|-\ell(w)}$ times the number of sequences $w =
w_1 \gtr w_2 \gtr \cdots \gtr w_m$ such that $w_m = w_\beta$ is a
Grassmannian permutation for $\beta$ and $w_i$ is not Grassmannian
for $i < m$.
\end{thm}

Let $w$ be any permutation and let $F = L_1 \oplus \dots \oplus L_f$
and $E = M_1 \oplus \cdots \oplus M_e$ be vector bundles on $\X$ which
are both direct sums of line bundles.  Buch
\cite{buch:littlewood-richardson} defines
\[ G_w(F-E) = G_w(1-L_1^{-1},\ldots, 1-L_{f}^{-1}; 
   1-M_1,\ldots,1-M_e)\in K(\X) \,.
\]
Since $G_w$ is symmetric, this definition extends to the case where $E$ and
$F$ do not split as direct sums.  Alternatively, using the identity
$\wedge^i(F^\vee) = (\wedge^{f-i} F)/(\wedge^f F)$, we may write $G_w(F-E)$ as
a Laurent polynomial in the exterior powers of $E$ and $F$, where only the top
power of $F$ is inverted.  As $G_w(1-e^{-X};1-e^{Y})$ is super-symmetric we
have $G_w(F\oplus H-E\oplus H)=G_w(F-E)$ for any third vector bundle $H$ on
$\X$.
% This allows one to regard $G_w$ as a well-defined function $G_w :
% K(\X) \to K(\X)$.
Finally, notice that
\begin{equation}
\label{combid}
G_{\alpha}(F-E)= \sum_{\beta,\gamma} d_{\beta\gamma}^{\alpha}\,
G_{\beta}(F-H)\, G_{\gamma}(H-E) \,,
\end{equation}
which follows from (\ref{refequ}) together with the super-symmetry
property.

%%%%%%%%%%%%%%%%%%%%%%%%%%%%%%%%%%%%%%%%%%%%%%%%%%%%%%%%%%%%%%%%%%%%%%
%%%   3 .  UNIVERSAL GROTHENDIECK POLYNOMIALS
%%%%%%%%%%%%%%%%%%%%%%%%%%%%%%%%%%%%%%%%%%%%%%%%%%%%%%%%%%%%%%%%%%%%%%

\section{Universal Grothendieck polynomials}
\label{UGP}

Consider a sequence
\[ 
B_\bull \ :\ B_1 \to B_2 \to \dots \to B_n
\] 
of vector bundles
and bundle maps over a non-singular variety $\X$.  Given rank
conditions $r = \{ r_{ij} \}$ for $1 \leq i < j \leq n$ there is a
{\em quiver variety\/} given by
\[ \Omega_r(B_\bull) = \{ x \in \X \mid 
   \rank(B_i(x) \to B_j(x)) \leq r_{ij} ~\forall i<j \} \,.
\]

For convenience, we set $r_{ii} = \rank (B_i)$ for all $i$, and we
demand that the rank conditions satisfy $r_{ij} \geq \max(r_{i-1,j},
r_{i,j+1})$ and $r_{ij} + r_{i-1,j+1} \geq r_{i-1,j} + r_{i,j+1}$ for
all $i \leq j$.  In this case, the expected codimension of
$\Omega_r(B_\bull)$ is the number $d(r) = \sum_{i<j}
(r_{i,j-1}-r_{ij})(r_{i+1,j}-r_{ij})$.  The main result of
\cite{buch:grothendieck} states that when the quiver variety
$\Omega_r(B_\bull)$ has this codimension, the class of its structure
sheaf is given by the formula

\begin{equation} \label{eqn:quiverformula}
  [{\mathcal O}_{\Omega_r(B_\bull)}] = \sum_\lambda c_\lambda(r)\, 
  G_{\lambda^1}(B_2-B_1) \cdots G_{\lambda^{n-1}}(B_n-B_{n-1}) \,.
\end{equation}
Here the sum is over finitely many sequences of partitions $\lambda =
(\lambda^1, \dots, \lambda^{n-1})$ such that $|\lambda| = \sum
|\lambda^i|$ is greater than or equal to $d(r)$.  The coefficients
$c_\lambda(r)$ are integers called {\em quiver coefficients}; they can
be computed by a combinatorial algorithm which we will not reproduce
here.  These coefficients are uniquely determined by the condition
that (\ref{eqn:quiverformula}) is true for all varieties $\X$ and
sequences $B_\bull$, as well as the condition that $c_\lambda(r) =
c_\lambda(r')$, where $r'=\{r'_{ij}\}$ is the set of rank conditions
given by $r'_{ij} = r_{ij}+1$ for all $i \leq j$.  Buch has
conjectured that the signs of these coefficients alternate with
codimension, that is, $(-1)^{|\lambda| - d(r)} c_\lambda(r) \,\geq\,
0$.

Suppose the index $p$ is such that all rank conditions $\rank(B_i(x)
\to B_p(x)) \leq r_{ip}$ and $\rank(B_p(x) \to B_j(x)) \leq r_{pj}$
may be deduced from other rank conditions.  As in 
\cite[\S 4]{buch.fulton:chern},
we will then say that the bundle $B_p$ is {\em inessential}.  Omitting
an inessential bundle $B_p$ from $B_\bull$ produces a sequence
\[
B'_\bull\ : \ B_1 \to \dots \to B_{p-1} \to B_{p+1} \to \dots \to B_n \,,
\]
where the map from $B_{p-1}$ to $B_{p+1}$ is the composition $B_{p-1}
\to B_p \to B_{p+1}$.  If $r'$ denotes the restriction of the rank
conditions to $B'_\bull$, we have that $\Omega_{r'}(B'_\bull) =
\Omega_r(B_\bull)$.  We can use (\ref{combid}) to expand any factor
$G_{\alpha}(B_{p+1}-B_{p-1})$ occurring in the quiver formula for
$\Omega_{r'}(B'_\bull)$ into a linear combination of products of the
form $G_{\beta}(B_p-B_{p-1}) G_{\gamma}(B_{p+1}-B_p)$, and thus arrive
at the quiver formula (\ref{eqn:quiverformula}) for
$\Omega_r(B_\bull)$.

The loci $\Omega_{w}(E_\bull\to F_\bull)$ of (\ref{eqn:univlocus}) are special
cases of these quiver varieties.  Given $w \in S_{n+1}$ we define rank
conditions $r^{(n)} = \{ r^{(n)}_{ij} \}$ for $1 \leq i \leq j \leq 2n$ by
\[ r^{(n)}_{ij} = \begin{cases}
   r_w(2n+1-j, i) & \text{if $i\leq n < j$} \\
   i & \text{if $j \leq n$} \\
   2n+1-j & \text{if $i \geq n+1$.}
\end{cases} \]
Then $\Omega_w(E_\bull \to F_\bull)$ is identical to the quiver
variety $\Omega_{r^{(n)}}(E_\bull \to F_\bull)$, and moreover we
have $d(r)=\ell(w)$.  We let $c^{(n)}_{w,\lambda} =
c_\lambda(r^{(n)})$ denote the quiver coefficients corresponding to
this locus.

Given vector bundles $E_1, \dots, E_n$ and $F_1, \dots, F_n$ on $\X$ we
define
\begin{equation} \label{eqn:quiv2univ}
  \G^{(n)}_w(F_\bull; E_\bull) = \sum_\lambda c_{w,\lambda}^{(n)}\, 
  G_{\lambda^1}(E_2-E_1) \cdots G_{\lambda^n}(F_n-E_n) \cdots 
  G_{\lambda^{2n-1}}(F_1-F_2) \,.
\end{equation}
It follows that $[\O_{\Omega_w}] = \G^{(n)}_w(F_\bull; E_\bull)$ when
the bundles are part of a sequence (\ref{eqn:univseq}) and the
codimension of $\Omega_w$ is equal to $\ell(w)$.

By definition, $\G^{(n)}_w(F_\bull; E_\bull)$ is a Laurent polynomial in the
exterior powers of the bundles $E_i$ and $F_i$, where only the top powers are
inverted.  We will call these polynomials {\em universal Grothendieck
polynomials}, in analogy with the term `universal Schubert polynomials' which
Fulton \cite{fulton:usp} used for his cohomology formula for $\Omega_w$.
The next lemma shows that the polynomial $\G_w(F_\bull; E_\bull)$ is
independent of $n$.  We will therefore drop this letter from the notation and
write simply $\G_w(F_\bull; E_\bull) = \G_w^{(n)}(F_\bull; E_\bull)$ when $w
\in S_{n+1}$.

\begin{lemma} \label{lemma:indep_n}
  Let $w \in S_{n+1}$.  The polynomial $\G_w^{(n+1)}(F_\bull; E_\bull)$
  is independent of $E_{n+1}$ and $F_{n+1}$ and agrees with
  $\G_w^{(n)}(F_\bull; E_\bull)$.
\end{lemma}
\begin{proof}
  Let $\X$ be a non-singular variety with a bundle sequence 
\[E_1 \to \dots \to E_{n+1} \to F_{n+1} \to \dots \to F_1
\]
such that the degeneracy locus $\Omega_w$ determined by this sequence has the
  expected codimension.  Since the same degeneracy locus is obtained by using
  the subsequence which skips the two middle bundles $E_{n+1}$ and $F_{n+1}$,
  it follows that $\G_w^{(n+1)}(F_\bull; E_\bull) = [\Omega_w] =
  \G_w^{(n)}(F_\bull; E_\bull)$, so the polynomials agree when evaluated in
  the Grothendieck ring $K(\X)$.
  
  To obtain the identity of polynomials, we need to construct a
  variety $\X$ such that all monomials in exterior powers which occur
  in either polynomials are linearly independent in $K(\X)$.
  Here we can use that on a Grassmannian $\Gr(m,N)$, all monomials of
  total degree at most $N/m-1$ in the exterior powers of the
  tautological subbundle are linearly independent.  Therefore we can
  take a product of Grassmannians
\[ {\mathcal Z} = \Gr(1,N) \times \dots \times \Gr(n+1,N) \times \Gr(n+1,N)
   \times \dots \times \Gr(1,N) \,,
\]
and let $\X$ be the bundle $\Hom(E_1,E_2) \oplus \dots \oplus
\Hom(E_{n+1},F_{n+1}) \oplus \dots \oplus \Hom(F_2,F_1)$, where the
bundles $E_i$ and $F_i$ are the tautological subbundles on ${\mathcal Z}$.  
When $N$ is sufficiently large, this variety $\X$ fits our purpose.
\end{proof}

In the remainder of this paper we will use without comment that the
universal Grothendieck polynomial $\G_w(F_\bull; E_\bull)$ is
determined by its values, as in the above proof.

%%%%%%%%%%%%%%%%%%%%%%%%%%%%%%%%%%%%%%%%%%%%%%%%%%%%%%%%%%%%%%%%%%%%%%
%%%   4 .  PROOF OF THE CAUCHY IDENTITY
%%%%%%%%%%%%%%%%%%%%%%%%%%%%%%%%%%%%%%%%%%%%%%%%%%%%%%%%%%%%%%%%%%%%%%

\section{Proof of the Cauchy identity}
\label{CF}

In this section we prove the Cauchy identity for universal
Grothendieck polynomials (Theorem \ref{thm:main2}).  We will assume
that $\X$ is a non-singular variety equipped with vector bundles $E_i$
and $F_i$ for $i \geq 1$, with $\rank E_i = \rank F_i = i$.

\begin{prop} \label{prop:push}
  Let $\pi : \Y = \Fl(E_n) \to \X$ be the bundle of flags in $E_n$,
  with tautological flag $0\subset U_1 \subset U_2 \subset \dots
  \subset U_n = \pi^*(E_n)$, and set
\[ \O_\Z = \prod_{1 \leq i \leq n-1} G_{(1^i)}(U_{i+1}/U_i - E_i)
   ~\in~ K(\Y) \,.
\]
Then we have 
\[ \pi_*(\G_w(F_\bull ; U_\bull) \cdot \O_Z) = \G_w(F_\bull ; E_\bull)
   ~\in~ K(\X) \,.
\]
\end{prop}
\begin{proof}
  Set $\wt \X = \Hom(E_1,E_2) \oplus \dots \oplus \Hom(E_n,F_n)
  \oplus \dots \oplus \Hom(F_2,F_1)$ and $\wt \Y = \Y \times_\X \wt{\X}$.
  It is enough to prove the proposition for the projection $\rho :
  \wt \Y \to \wt \X$.  Notice that on $\wt \Y$ we have a
  universal bundle sequence $E_\bull \to F_\bull$, as well as the
  tautological flag $U_\bull \subset E_n$.
  
  Let $\Z_{n-1} = Z(E_{n-1} \to U_n/U_{n-1}) \subset \wt \Y$.  On this
  locus the map $E_{n-1} \to U_n$ factors through $U_{n-1}$.  We then
  set $\Z_{n-2} = Z(E_{n-2} \to U_{n-1}/U_{n-2}) \subset \Z_{n-1}$ and
  inductively $\Z_i = Z(E_i \to U_{i+1}/U_i)$ for $i = n-1,\dots,2,1$.
  Notice that the structure sheaf of $\Z = \Z_1$ is given by the
  expression of the proposition (see e.g.\ 
  \cite[Thm.~2.3]{buch:grothendieck}).  Now $\rho$ maps the locus
  $\Omega_w(U_\bull \to F_\bull) \cap \Z \subset \wt \Y$ birationally
  onto $\Omega_w(E_\bull \to F_\bull) \subset \wt \X$.
  In fact, the open subset of $\Z$ where each map $E_i \to U_i$ is an
  isomorphism maps isomorphically to the open subset of $\wt \X$ where
  all maps $E_{i-1} \to E_i$ are bundle inclusions, and furthermore
  these subsets meet the given (irreducible) degeneracy loci in $\Z$
  and $\wt \X$.  This implies the desired result, because all involved
  degeneracy loci are Cohen-Macaulay with rational singularities
  \cite{lakshmibai.magyar:degeneracy} and have their expected
  codimensions.
\end{proof}

This proposition allows us to prove a special case of the Cauchy formula,
arguing as in \cite[\S 3]{fulton:usp}.  We let $\C^\bull$ denote a
sequence of trivial bundles.  When used in a polynomial, the exterior power
$\wedge^i \C^m$ equals the binomial coefficient $\binom{m}{i}$.

\begin{cor} \label{cor:specialcauchy}
We have
\[ \G_w(F_\bull; E_\bull) = \sum_{u\cdot v = w} (-1)^{\ell(uvw)}\,
   \G_u(\C^\bull; E_\bull)\, \G_v(F_\bull; \C^\bull) \,.
\]
\end{cor}
\begin{proof}
  Let $\pi : \Y = \Fl(E_n) \to \X$ be the bundle of flags in $E_n$, with
  tautological flag $U_\bull$ as in Proposition \ref{prop:push}.  Assume at
  first that the bundles $F_\bull$ form a (descending) complete flag.  Then by
  \cite[Thm.~2.1]{buch:grothendieck} (which generalizes
  \cite[Thm.~3]{fulton.lascoux:pieri}), we have
\[ \G_w(F_\bull;U_\bull) = 
   \Groth_w(1-L_1^{-1},\dots,1-L_n^{-1} \,;\, 1-M_1,\dots,1-M_n)
\]
in $K(\Y)$, where $L_i = \ker(F_i \to F_{i-1})$ and $M_i =
U_i/U_{i-1}$.  The Cauchy identity for double Grothendieck polynomials
(\ref{eqn:cauchy_old}) therefore implies that
\[ \G_w(F_\bull ; U_\bull) = \sum_{u\cdot v=w} (-1)^{\ell(uvw)}\,
   \G_u(\C^\bull ; U_\bull) \cdot \G_v(F_\bull ; \C^\bull) \,.
\]
By multiplying this identity by the class $\O_\Z$ of
\refprop{prop:push}, and pushing the result down to $\X$, we get the
identity of the theorem.

Now assume that the bundles $F_i$ are arbitrary.  By the case just
proved we have
\[\begin{split}
  \G_w(F_\bull ; U_\bull) 
  &= \G_{w^{-1}}(U_\bull^\vee ; F_\bull^\vee) \\
  &= \sum_{u\cdot v=w} (-1)^{\ell(uvw)}\, 
     \G_{v^{-1}}(\C^\bull;F_\bull^\vee) \cdot
     \G_{u^{-1}}(U_\bull^\vee;\C^\bull) \\
  &= \sum_{u\cdot v=w} (-1)^{\ell(uvw)}\, \G_u(\C^\bull;U_\bull)\cdot
     \G_v(F_\bull ; \C^\bull)
\end{split}\]
in $K(\Y)$.  After multiplying with $\O_\Z$, this identity
pushes forward to give the corollary in full generality.
\end{proof}

For the general case of the Cauchy formula we need the following
vanishing theorem.  We let $H_\bull$ denote a third collection of
vector bundles on $\X$, $\rank H_i = i$.

\begin{prop}
\label{prop:vanish}
Choose $m\geq 0$ and substitute $H_{j}$ for $F_{j}$ and $E_{j}$ in
$\G_{w}(F_\bull;E_\bull)$ for all $j\geq m+1$. We then have
\[
  \G_{w}(F_1,\ldots,F_m,H_{m+1},\ldots;E_1,\ldots,E_m,H_{m+1},\ldots)
  = \begin{cases} \G_w(F_\bull;E_\bull) & \text{if $w \in S_{m+1}$} \\
  0 & \text{otherwise.} \end{cases}
\]
\end{prop}
\begin{proof}
  If $w\in S_{m+1}$, then $\G_{w}(F_\bull,E_\bull)$ is independent of
  the bundles $F_{j}$ and $E_{j}$ for $j\geq m+1$ by
  Lemma~\ref{lemma:indep_n}.
  
  Assume that $w \in S_{n+1} \smallsetminus S_{n}$ where $n>m$.  We
  claim $\G_{w}(F_\bull,E_\bull)$ vanishes as soon as we set $F_n =
  E_n$.  To see this, let $\X$ be a variety with bundles $F_j$ for $1
  \leq j \leq n$ and $E_j$ for $1 \leq j \leq n-1$ such that all
  monomials in the polynomial $\G_w(F_\bull; E_1,\dots,E_{n-1},F_n)$
  are linearly independent.

  If we set $E_n = F_n$ then we have a sequence of bundles $E_\bull
  \to F_\bull$ for which the map $E_{n}\to F_{n}$ is the identity and
  all other maps are zero.  Since $r_{w}(n,n) = n-1$ it follows that
  the locus $\Omega_{w}(E_\bull \to F_\bull)$ is empty.  Since
  $\G_{w}(F_\bull;E_\bull)$ represents the class of the structure sheaf
  of this locus, it must be equal to zero.
\end{proof}

For any commutative ring $R$, let $R(S_{\infty})$ denote the
$R$-module of all functions on $S_{\infty} = \bigcup_n S_n$ with
values in $R$.  For $f, g \in R(S_\infty)$ we define the product
\begin{equation}
\label{eqn:product_RS}
(fg)(w) = \sum_{u\cdot v = w} (-1)^{\ell(uvw)} f(u) g(v)
\end{equation}
where (as always) the sum is over factorizations of $w$ in the
degenerate Hecke algebra.  It is straightforward to check that this
multiplication is associative and that the identity element is the
characteristic function ${\underbar{1}}$ of the identity permutation.
We will need the following variation of
\cite[(6.6)]{Macdonald:schub_notes}.

\begin{lemma} \label{lemma:Macdonald_modified}
  Let $f,g,h \in R(S_{\infty})$.  Assume that for any permutation $w
  \in S_\infty$, the sum $\sum_{u\cdot w=w} (-1)^{\ell(u)} f(u)$ is
  not a zero divisor in $R$.
\begin{itemize}
\item[(i)] If $fg=f$ then $g={\underbar{1}}$.
\item[(ii)] If $fh={\underbar 1}$ then $hf={\underbar 1}$.
\end{itemize}
\end{lemma}
\begin{proof} Since $f(1)g(1) = fg(1) = f(1)$ and $f(1) = \sum_{u\cdot
    1=1} (-1)^{\ell(u)} f(u)$ is a non-zero divisor, it follows that
  $g(1) = 1$.  Let $w \neq 1 \in S_\infty$ be given and assume
  inductively that $g(v) = 0$ for $0 < \ell(v) < \ell(w)$.  Notice
  that if $u \cdot v = w$ in the degenerate Hecke algebra then
  $\ell(v) \leq \ell(w)$, and this inequality is sharp if $v \neq w$.
  We therefore have $f(w) = fg(w) = f(w) + \left(\sum_{u\cdot w=w}
  (-1)^{\ell(u)} f(u)\right) g(w)$, which implies that $g(w) = 0$.  This
  proves (i), and (ii) follows by setting $g = hf$.
\end{proof}

% \begin{proof} First observe that
% \begin{equation} \label{eqn:easy_inequality}
%   \ell(u\cdot v)\geq {\rm max}(\ell(u),\ell(v)) \,.
% \end{equation}
% Thus if $u,v\in S_{\infty}$ and $u\cdot v=1$, then $u=v=1$. Hence
% $fg=f$ implies that $f(1)g(1)=f(1)$ and $g(1)=1$.
%
% We prove by induction on $\ell(w)$ that $g(w)=0$ for all $w\neq 1$.
% Let $r>0$ and assume $g(v)=0$ for all $v\in S_{\infty}$ such that
% $1\leq \ell(v)\leq r-1$. Let $w\in S_{\infty}$ be of length $r$. Then
% \[ f(w) = fg(w) = 
%    f(w)g(1)+f(1)g(w)+\sum_{u,v}(-1)^{\ell(uvw)}f(u)g(v) \,, 
% \]
% where the sum on the right is over $u,v\in S_{\infty}$ such that
% $u,v\neq 1$ and $u\cdot v=w$. By (\ref{eqn:easy_inequality}), $1\leq
% l(u)\leq r-1$, so $g(v)=0$. Therefore $f(1)g(w)=0$ and $g(w)=0$. Hence
% (i) holds.
%
% To prove (ii), notice from (\ref{eqn:easy_inequality}) we have
% $f(1)g(1)=1$ so $f(1)$ is a unit in $R$. Since $f(gf)=(fg)f=f$ we
% conclude $gf={\underline 1}$ by (i).
% \end{proof}

\begin{thmII}[Cauchy formula]
Let $E_{i}, F_{i}$, and $H_{i}$ for $i=1,\ldots, n$ be three collections
of vector bundles on $\X$.  Then for any $w \in S_{n+1}$ we have
\[ \G_w(F_\bull; E_\bull) = \sum_{u \cdot v=w} (-1)^{\ell(uvw)}\,
   \G_u(H_\bull; E_\bull)\, \G_v(F_\bull; H_\bull) \,.
\]
\end{thmII}
\begin{proof}
  Let $\G(F_\bull;E_\bull)$ denote the function from permutations to
  $K(\X)$ which maps $w$ to $\G_w(F_\bull;E_\bull)$.  Using the product
  (\ref{eqn:product_RS}) we have by Corollary~\ref{cor:specialcauchy}
  that
\[ \G(F_\bull; E_\bull) = 
   \G(\C^\bull; E_\bull) \G(F_\bull; \C^\bull) \,.
\]
Proposition~\ref{prop:vanish} implies that $\G(\C^\bull;H_\bull)
\G(H_\bull,\C^\bull)={\underline 1}$, and since $\G_w(\C^\bull;H_\bull)$
lies in the augmentation ideal of $K(\X)$ for $w \neq 1$, the function
$f = \G(\C^\bull;H_\bull)$ satisfies the requirement in
Lemma~\ref{lemma:Macdonald_modified}.  It follows that
$\G(H_\bull,\C^\bull)\G(\C^\bull,H_\bull)={\underline 1}$.  We conclude
that
\[ \G(F_\bull;E_\bull) 
   = \G(\C^\bull;E_\bull) \G(F_\bull;\C^\bull)
   = \G(\C^\bull;E_\bull) \G(H_\bull;\C^\bull) \G(\C^\bull;H_\bull)
     \G(F_\bull;\C^\bull) \,,
\]
and therefore $\G(F_\bull;E_\bull) = \G(H_\bull;E_\bull)
\G(F_\bull;H_\bull)$, as required.
\end{proof}

For later use, we notice that Theorem~\ref{thm:main2} and
Proposition~\ref{prop:vanish} together imply that
% begin{equation} \label{eqn:splitright}
%   P_w(F_\bull; E_\bull) = \sum_{u \cdot v = w} (-1)^{\ell(uvw)}
%   P_u(\C^1,\dots,\C^r,F_{r+1}, \dots ; E_\bull)
%   P_v(F_1,\dots,F_r; \C^\bull)
% \end{equation}
% and similarly,
% \begin{equation} \label{eqn:splitleft}
%   P_w(F_\bull; E_\bull) = \sum_{u \cdot v = w} (-1)^{\ell(uvw)}
%   P_u(\C^\bull; E_1,\dots,E_r)
%   P_v(F_\bull; \C^1, \dots, \C^r, E_{r+1}, \dots) \,.
% \end{equation}
\begin{equation} \label{eqn:split}
\begin{split}
&\G_w(F_\bull; E_\bull) \\
&= \sum_{u \cdot v = w} (-1)^{\ell(uvw)}
  \G_u(\C^1,\dots,\C^r,F_{r+1}, F_{r+2}, \dots \,;\, E_\bull)
  \G_v(F_1,\dots,F_r \,;\, \C^\bull) \\
&= \sum_{u \cdot v = w} (-1)^{\ell(uvw)}
  \G_u(\C^\bull \,;\, E_1,\dots,E_r)
  \G_v(F_\bull \,;\, \C^1, \dots, \C^r, E_{r+1}, E_{r+2}, \dots) \,.
\end{split}
\end{equation}

%%%%%%%%%%%%%%%%%%%%%%%%%%%%%%%%%%%%%%%%%%%%%%%%%%%%%%%%%%%%%%%%%%%%%%
%%%   5 .  MAIN THEOREM
%%%%%%%%%%%%%%%%%%%%%%%%%%%%%%%%%%%%%%%%%%%%%%%%%%%%%%%%%%%%%%%%%%%%%%

\section{Proof of Theorem~\ref{thm:main1}} 

In this section we derive Theorem~\ref{thm:main1} from the Cauchy
identity by using a $K$-theoretic version of the arguments found in
\cite{schquiv}.  In what follows, it will be
convenient to work with the element $P^{(n)}_w \in \Gamma^{\otimes
  2n-1}$ defined by
\[ P^{(n)}_w = \sum_{\lambda} c_{w,\lambda}^{(n)}\, 
   G_{\lambda^1} \otimes \dots \otimes G_{\lambda^{2n-1}} \,.
\]
With this notation, we can restate Theorem~\ref{thm:main1} as follows:

\begin{thm1a}
%\label{thm:oneprime}
  For any permutation $w \in S_{n+1}$ we have
\[ P_w^{(n)} = \sum_{u_1 \cdots u_{2n-1} = w} 
   (-1)^{\ell(u_1 \cdots u_{2n-1} w)}\,
   G_{u_1} \otimes \dots \otimes G_{u_{2n-1}} 
\]
in $\Gamma^{\otimes 2n-1}$, where the sum is over all factorizations $w = u_1
\cdots u_{2n-1}$ in the degenerate Hecke algebra such that $u_i \in
S_{\min(i,2n-i)+1}$ for each $i$.
\end{thm1a}
\begin{proof}
  Since $r_w(p,q)+m = r_{1^m\times w}(p+m,q+m)$ for $m \geq 0$, it
  follows that the coefficients $c^{(n)}_{w,\lambda}$ are uniquely
  defined by the condition that
\begin{multline} \label{eqn:shift2quiv}
  \G_{1^m\times w}(F_\bull;E_\bull) = \sum_\lambda
  c_{w,\lambda}^{(n)}\, G_{\lambda^1}(E_{2+m}-E_{1+m}) \cdots \\ 
  G_{\lambda^n}(F_{n+m}-E_{n+m}) \cdots
  G_{\lambda^{2n-1}}(F_{1+m}-F_{2+m})
\end{multline}
for all $m \geq 0$ (see \cite{buch:stanley} and also the discussion
after the proof of Theorem~4.1 in \cite{buch:grothendieck}).

Given any two integers $p \leq q$ we let $P_w^{(n)}[p,q]$ denote the
sum of the terms of $P^{(n)}_w$ for which $\lambda^i$ is empty when
$i<p$ or $i>q$:
\begin{equation*}
  P^{(n)}_w[p,q] = \sum_{\lambda: \lambda^i = \emptyset 
  \text{ for } i \not \in [p,q]} 
  c_{w,\lambda}^{(n)}\, G_{\lambda^1} \otimes \dots \otimes 
  G_{\lambda^{2n-1}} \,.
\end{equation*}

\begin{lemma} \label{lemma:splitquiv}
For any $1 < i \leq 2n-1$ we have
\[ P^{(n)}_w = 
   \sum_{u \cdot v = w} (-1)^{\ell(uvw)} 
   P^{(n)}_u[1,i-1]\cdot P^{(n)}_v[i,2n-1] \,. 
\]
\end{lemma}
\begin{proof}
  We will do the case $i \leq n$, the other one is similar.  For any
  element $f = \sum c_\lambda\, G_{\lambda^1} \otimes \dots \otimes
  G_{\lambda^{2N-1}} \in \Gamma^{\otimes 2N-1}$, we set
\[ f(F_\bull;E_\bull) = \sum c_\lambda\, G_{\lambda^1}(E_2-E_1) \cdots
   G_{\lambda^N}(F_N-E_N) \cdots G_{\lambda^{2N-1}}(F_1-F_2) \,.
\]
Equation (\ref{eqn:shift2quiv}) implies that $P^{(n)}_w \in
\Gamma^{\otimes 2n-1}$ is the unique element satisfying that
$(1^{\otimes m} \otimes P^{(n)}_w \otimes 1^{\otimes m})
(F_\bull;E_\bull) = \G_{1^m\times w}(F_\bull;E_\bull)$ in $K(\X)$ for
all $m$.  This uniqueness is preserved even if we make $E_{i+m}$
trivial.  The right hand side of the identity of the lemma satisfies
this by equation (\ref{eqn:split}) applied to $1^m \times w$.
\end{proof}

\begin{lemma} \label{lemma:stable}
For $1 \leq i \leq 2n-1$ we have
\[ P_w^{(n)}[i,i] = \begin{cases}
   1^{\otimes i-1} \otimes G_w \otimes 1^{\otimes 2n-1-i}
   & \text{if $w \in S_{m+1}$, $m=\min(i,2n-i)$}, \\
   0 & \text{otherwise.}
\end{cases} \]
\end{lemma}
\begin{proof}
  For simplicity we will assume that $m = i$.  If $w \not \in S_{m+1}$
  then it follows from Proposition~\ref{prop:vanish} or the algorithm
  for quiver coefficients of \cite[\S 4]{buch:grothendieck} that
  $P_w^{(n)}[1,m] = 0$, which proves the lemma.  Assume now that $w
  \in S_{m+1}$.  It is proved in \cite[(5.2)]{buch:grothendieck} that
  $P_w^{(m)}[m,m] = 1^{\otimes m-1} \otimes G_w \otimes 1^{\otimes
    m-1}$.  Let $\Phi : \Gamma^{\otimes 2m-1} \to \Gamma^{\otimes
    2n-1}$ be the linear map given by
\[ \Phi(G_{\lambda^1}\otimes\dots\otimes G_{\lambda^{2m-1}}) =
   G_{\lambda^1} \otimes \dots \otimes G_{\lambda^{m-1}} \otimes
   \Delta^{2n-2m}(G_{\lambda^m}) \otimes G_{\lambda^{m+1}} \otimes 
   \dots \otimes G_{\lambda^{2m-1}}
\]
where $\Delta^{2n-2m} : \Gamma \to \Gamma^{\otimes 2n-2m+1}$ denotes
the $(2n-2m)$-fold coproduct, that is,
\[ \Delta^{2n-2m}(G_{\lambda^m}) =
   \sum_{\tau_1,\dots,\tau_{2n-2m+1}}
   d^{\lambda^m}_{\tau_1,\dots,\tau_{2n-2m+1}}\,
   G_{\tau_1} \otimes \dots \otimes G_{\tau_{2n-2m+1}}
\]
(see \cite[Corollary~6.10]{buch:littlewood-richardson}). In the
definition of the locus $\Omega_w(E_\bull \to F_\bull)$, the bundles
$F_i$ and $E_i$ for $i \geq m+1$ are inessential in the sense of
Section \ref{UGP}, which implies that $\Phi(P^{(m)}_w) =
P^{(n)}_w$.  Now the result follows from the identity
$P^{(n)}_w[m,2n-m] = \Phi(P^{(m)}_w[m,m]) = 1^{\otimes m-1} \otimes
\Delta^{2n-2m}(G_w) \otimes 1^{\otimes m-1}$.
\end{proof}

Theorem $1'$ follows immediately from
Lemma~\ref{lemma:splitquiv} and Lemma~\ref{lemma:stable}.
\end{proof}

\begin{cor} \label{cor:quivcoef}
\label{altcor}
  Let $w\in S_{n+1}$ and let ${\lambda} = (\lambda^1, \lambda^2, \ldots,
  \lambda^{2n-1})$ be a sequence of partitions.  Then we have
\[ c_{w,\lambda}^{(n)} = (-1)^{|\lambda|-\ell(w)}
  \sum_{u_1 \cdots u_{2n-1} = w} 
  |\, a_{u_1,\lambda^1}\, a_{u_2,\lambda^2} \cdots
  a_{u_{2n-1},\lambda^{2n-1}} |
\]
where $|\lambda|=\sum |\lambda^i|$, $a_{u_i,\lambda^i}$ is the
coefficient of $G_{\lambda^i}$ in $G_{u_i} \in \Gamma$, and the sum is
over all factorizations of $w$ in the degenerate Hecke algebra such
that $u_i \in S_{\min(i,2n-i)+1}$ for each $i$.
\end{cor}

Since Lascoux's formula (Theorem~\ref{thm:lascoux}) implies that
$a_{u_i,\lambda^i} = (-1)^{|\lambda^i|-\ell(u_i)}\, |a_{u_i,\lambda^i}|$, 
Corollary \ref{altcor} follows immediately from Theorem~\ref{thm:main1}.  
This verifies the
alternation of signs for the quiver coefficients $c_{w,\lambda}^{(n)}$, which
was conjectured in \cite{buch:grothendieck}.  In addition, by combining
the above corollary with Lascoux's formula we obtain an explicit
combinatorial formula for these coefficients.

%%%%%%%%%%%%%%%%%%%%%%%%%%%%%%%%%%%%%%%%%%%%%%%%%%%%%%%%%%%%%%%%%%%%%%
%%%   1 .  SPLITTING GROTHENDIECK POLYNOMIALS
%%%%%%%%%%%%%%%%%%%%%%%%%%%%%%%%%%%%%%%%%%%%%%%%%%%%%%%%%%%%%%%%%%%%%%

\section{Splitting Grothendieck polynomials}
\label{SGP}

In this section we specialize universal Grothendieck polynomials to
the double Grothendieck polynomials of Lascoux and Sch{\"u}tzenberger.
This leads to new expressions for double Grothendieck polynomials in
terms of quiver coefficients.  Recall that a permutation $w$ has a
descent at position $i$ if $w(i)>w(i+1)$. We say that a sequence
$\{a_k\}:a_1<\ldots<a_p$ of integers is {\em compatible with} $w$ if
all descent positions of $w$ are contained in $\{a_k\}$.
%Theorem~\ref{thm:main3} is immediate by taking $a_i=b_i =
%i$ for $i=1,2,\ldots$ in part (II) of the following theorem:

\begin{thm}
\label{thm:main4}
Let $w\in S_{n+1}$ and let $1 \leq a_1 < \dots < a_p \leq n$ and $1
\leq b_1 < \dots < b_q \leq n$ be two sequences compatible with $w$
and $w^{-1}$, respectively, and set $X_i = \{x_{a_{i-1}+1}, \dots,
x_{a_i}\}$ and $Y_i = \{y_{b_{i-1}+1}, \dots, y_{b_i}\}$. Then we
have
\begin{equation} \label{eqn:desired}
  \Groth_w(X;Y)=\sum_\mu {\wt c}_{w,\mu}\, G_{\mu^1}(X_p;0) \cdots 
  G_{\mu^p}(X_1;Y_1) \cdots G_{\mu^{p+q-1}}(0;Y_q) \,,
\end{equation}
where the sum is over sequences of partitions $\mu = (\mu^1, \ldots,
\mu^{p+q-1})$, and ${\wt c}_{w,\mu}$ is the quiver coefficient
$c^{(n)}_{w_0 w^{-1} w_0,\lambda}$, where $w_0\in S_{n+1}$ is the
longest permutation and $\lambda = (\lambda^1, \ldots,
\lambda^{2n-1})$ is given by
\[ \lambda^i = \begin{cases} \mu^k & \text{if $i = a_{k+1}-1$} \\
  \mu^p & \text{if $i = n$} \\
  \mu^{p+q-k} & \text{if $i = 2n-b_{k+1}+1$} \\
  \emptyset & \text{otherwise.}
  \end{cases}
\]
\end{thm}
\begin{proof}
Let $V$ be a vector bundle of rank $n+1$ and let 
\[ E_1 \subset E_2 \subset \dots \subset E_n \subset V 
   \tto F_n \tto \cdots \tto F_2 \tto F_1
\]
be a complete flag followed by a dual complete flag of $V$.  By
\cite[Thm.~2.1]{buch:grothendieck}, the class of the structure sheaf
of $\Omega_w(E_\bull \to F_\bull)$ is given by $\Groth_w(X;Y)$, where
we set $x_i = 1 - [\ker(F_i \to F_{i-1})]^{-1}$ and $y_i = 1 -
[E_i/E_{i-1}]$ in $K(\X)$.

Set $E_{i}'=V/E_i$ and $F_{i}'=\ker(V\tto F_i)$. This yields the
sequence
\[ F'_n \subset \dots \subset F'_1 \subset V  
   \tto E'_1 \tto \cdots \tto E'_n
\]
and it is easy to check that $\Omega_w(E_\bull\to F_\bull)=\Omega_{w_0
  w^{-1} w_0}(F_\bull'\to E_\bull')$ as subschemes of $\X$, where
$w_0$ is the longest permutation in $S_{n+1}$.

Define a third bundle sequence $\wt E'_\bull \to \wt F'_\bull$ as
follows.  For $a_{k-1} < i \leq a_k$ we set $\wt F'_i = F'_{a_k}
\oplus \C^{a_k-i}$ and for $b_{k-1} < i \leq b_k$ we set $\wt E'_i =
E'_{b_k} \oplus \C^{b_k-i}$.  The maps of the sequence $\wt E'_\bull
\to \wt F'_\bull$ can be chosen arbitrarily so that the subsequence
\[ 
\wt F'_{a_p} \to \dots \to \wt F'_{a_1} \to \wt E'_{b_1} \to \dots
    \to \wt E'_{b_q}
\]
agrees with the corresponding subsequence of $F'_\bull \to E'_\bull$,
the map $\wt F'_{i+1} \to \wt F'_i$ is an inclusion of vector bundles
for $i \not \in \{a_k\}$, and $\wt E'_i \to \wt E'_{i+1}$ is
surjective for $i \not \in \{b_k\}$.  Now \cite[Lemma 3]{schquiv}
implies that $\Omega_{w_0 w^{-1} w_0}(F'_\bull \to E'_\bull) =
\Omega_{w_0 w^{-1} w_0}(\wt F'_\bull \to \wt E'_\bull)$.  These
identities of schemes show that
\[ 
\Groth_w(X;Y) = \G_{w_0 w^{-1} w_0}(\wt E'_\bull; \wt F'_\bull) \,. 
\]
Equation (\ref{eqn:desired}) now follows from equation
(\ref{eqn:quiv2univ}).  In fact, $G_\alpha(\wt F'_i - \wt F'_{i+1})$
is non-zero only if $\alpha$ is empty or $i = a_k$ for some $k$, and
when $i = a_k$ we have $G_\alpha(\wt F'_i - \wt F'_{i+1}) =
G_\alpha(X_{k+1})$.  Similarly, $G_\alpha(\wt E'_{i+1} - \wt E'_i)$ is
zero unless $\alpha$ is empty or $i = b_k$ for some $k$, and for $i =
b_k$ we have $G_\alpha(\wt E'_{i+1} - \wt E'_i) =
G_\alpha(0;Y_{k+1})$.  Finally $G_\alpha(\wt E'_n - \wt F'_n) =
G_\alpha(X_1;Y_1)$.

This proves (\ref{eqn:desired}) in the Grothendieck ring $K(\X)$, in
which there are relations between the variables $x_i$ and $y_i$
(including e.g.\ the relations $e_j(x_1,\dots,x_{n+1}) =
e_j(y_1,\dots,y_{n+1})$ for $1 \leq j \leq n+1$).  We claim, however,
that (\ref{eqn:desired}) holds as an identity of polynomials in
independent variables.  For this, one checks that the definition of
$\wt c_{w,\mu}$ is independent of $n$, i.e.\ the coefficient
$c^{(n)}_{w_0 w^{-1} w_0,\lambda}$ does not change when $n$ is
replaced with $n+1$ and $w_0$ with the longest element in $S_{n+2}$.
If we choose $n$ sufficiently large, we can construct a variety $\X$
on which (\ref{eqn:desired}) is true, and where all relevant monomials
in the variables $x_i$ and $y_i$ are linearly independent.  This
establishes the claim.
\end{proof}

It follows from Theorem \ref{thm:main4} that the monomial coefficients
of Grothendieck polynomials are special cases of the $K$-theoretic
quiver coefficients $c^{(n)}_{w,\lambda}$.  Explicit formulas for the
monomial coefficients in terms of `non-reduced RC-graphs' are one of
the many consequences of Fomin and Kirillov's work
\cite{fomin.kirillov:grothendieck, fomin.kirillov:yang-baxter}.  We
will finish this paper by proving a different formula which
generalizes \cite[Thm.~1.1]{BJS} and \cite[Cor.~4]{schquiv}.

\begin{lemma} \label{lemma:onevar}
  Let $w$ be a permutation and $p \geq 0$ an integer.  Then the
  coefficient $a_{w,(p)}$ of Theorem \ref{thm:lascoux} is given by
\[ a_{w,(p)} = \begin{cases} 1 & \text{if $w = s_{i_1} \cdots s_{i_p}$
  for integers $i_1 > \dots > i_p$}, \\
  0 & \text{otherwise.} \end{cases}
\]
\end{lemma}
\begin{proof}
  Let $x$ be a variable and consider the degenerate Hecke algebra
  tensored with ${\mathbb Z}[x]$.  It follows from
  \cite[Thm.~2.3]{fomin.kirillov:grothendieck} that the Grothendieck
  polynomial $\Groth_w(x) = \Groth_w(x,0,\dots \,;\, 0,0,\dots)$ is
  equal to the coefficient of $w$ in the expansion of the product
\[ (1+x s_n)(1+x s_{n-1}) \cdots (1+x s_1) \]
in this algebra.  In other words, $\Groth_w(x)$ is non-zero exactly
when $w$ has a decreasing reduced word, in which case we have
$\Groth_w(x) = x^{\ell(w)}$.  The same is therefore true for the
stable polynomial $G_w(x)$.  The lemma follows from this because
$G_\beta(x) = 0$ for any partition $\beta$ of length at least two,
while $G_{(p)}(x) = x^p$.
\end{proof}

Using Fomin's identity $G_w(X;Y) = G_{w_0 w w_0}(Y;X)$ (see
\cite[Lemma~3.4]{buch:littlewood-richardson}) we similarly obtain that
$a_{w,(1^p)} = a_{w_0 w w_0,(p)}$ is equal to one if $w$ has an
increasing reduced word of length $p$, and $a_{w,(1^p)} = 0$
otherwise.

\begin{cor} \label{cor:groth_bjs}
  Let $w \in S_n$, let $x^u y^v = x_1^{u_1} \cdots x_{n-1}^{u_{n-1}}
  y_1^{v_1} \cdots y_{n-1}^{v_{n-1}}$ be a monomial, and set $g_i =
  \sum_{k=n-i}^{n-1} v_k$, $f_i = g_{n-1} + \sum_{k=1}^i u_k$, and
  $r=f_{n-1}=|u|+|v|$.  Then the coefficient of $x^u y^v$ in the
  double Grothendieck polynomial $\Groth_w(X;Y)$ is equal to
  $(-1)^{r-\ell(w)}$ times the number of factorizations $w=s_{e_1}
  \cdots s_{e_r}$ in the degenerate Hecke algebra such that $n-i\leq
  e_{g_{i-1}+1} < \dots < e_{g_i}$ and $e_{f_{i-1}+1} > \dots >
  e_{f_i} \geq i$ for all $1 \leq i \leq n-1$.
\end{cor}
\begin{proof}
  We apply Theorem \ref{thm:main4} to $\sigma = 1 \times w$ with
  $p=q=n$ and $a_i=b_i=i$, and use that $\Groth_w(X;Y) =
  \Groth_\sigma(0,x_1,\dots,x_{n-1} ; 0,y_1,\dots,y_{n-1})$.  The
  coefficient of $x_2^{u_1} \cdots x_n^{u_{n-1}} y_2^{v_1} \cdots
  y_n^{v_{n-1}}$ in $\Groth_\sigma(X;Y)$ is equal to $\tilde
  c_{\sigma,\lambda} = c_{w_0 \sigma^{-1} w_0,\lambda}$ where $\lambda
  = ((u_{n-1}), \dots, (u_1), \emptyset, (1^{v_1}), \dots,
  (1^{v_{n-1}}))$.  By Corollary~\ref{cor:quivcoef} and
  Lemma~\ref{lemma:onevar} this coefficient is equal to $\pm$ the
  number of factorizations $w_0 \sigma^{-1} w_0 = \tau_1 \cdots
  \tau_{n-1} \tau_{n+1} \cdots \tau_{2n-1}$ such that each $\tau_i$ is
  in $S_{\min(i,2n-i)+1}$ and has a decreasing reduced word of length
  $u_i$ for $i < n$ and an increasing reduced word of length
  $v_{2n-i}$ for $i > n$.  The sequences $(e_1,\dots,e_r)$ of the
  corollary are the corresponding factorizations of $w$.
\end{proof}

\begin{example} The double Grothendieck polynomials for the elements $s_i$
of length one in $S_n$ are given by the formula
\[
\Groth_{s_i}(X;Y) = \sum_{\delta}(-1)^{|\delta|-1}(xy)^{\delta}=
\sum_{\delta} (-1)^{|\delta|-1}x_1^{\delta_1}
\cdots x_{n-1}^{\delta_{n-1}}y_1^{\delta_n}\cdots
y_{n-1}^{\delta_{2n-2}}
\]
where the sum is over the $4^{n-1}-1$ strings $\delta=(\delta_1,\ldots,
\delta_{2n-2})$ with $\delta_i\in\{0,1\}$ for each $i$ and $|\delta|=
\sum\delta_i>0$. 
For instance, $\Groth_{s_1}(X;Y) = x_1+y_1-x_1y_1$ and
\begin{eqnarray*}
\Groth_{s_2}(X;Y) &=&
x_1+x_2+y_1+y_2-x_1x_2-x_1y_1-x_1y_2-x_2y_1-x_2y_2-y_1y_2 \\
&&+\  x_1x_2y_1+x_1x_2y_2+x_1y_1y_2+x_2y_1y_2-x_1x_2y_1y_2.
\end{eqnarray*}
This follows from Corollary \ref{cor:groth_bjs} since the
factorizations of $s_i$ in the degenerate Hecke algebra are exactly
the non-zero powers of $s_i$.
\end{example}

%%%%%%%%%%%%%%%%%%%%%%%%%%%%%%%%%%%%%%%%%%%%%%%%%%%%%%%%%%%%%%%%%%%%%%
%%%   BIBLIOGRAPHY
%%%%%%%%%%%%%%%%%%%%%%%%%%%%%%%%%%%%%%%%%%%%%%%%%%%%%%%%%%%%%%%%%%%%%%

\providecommand{\bysame}{\leavevmode\hbox to3em{\hrulefill}\thinspace}
\providecommand{\MR}{\relax\ifhmode\unskip\space\fi MR }
\providecommand{\MRhref}[2]{%
  \href{http://www.ams.org/mathscinet-getitem?mr=#1}{#2}
}
\providecommand{\href}[2]{#2}

%\bibliographystyle{amsplain}
%\bibliography{../BibTeX/database,ugp}

\begin{thebibliography}{99}

\bibitem{BJS} S. Billey, W. Jockusch and R. P. Stanley,
{\em Some combinatorial properties of Schubert polynomials},
J. Algebraic Combin. {\bf 2} (1993), no. 4, 345--374. 

\bibitem{buch:stanley}
A.~S. Buch, \emph{Stanley symmetric functions and quiver varieties,}
J.~of Algebra {\bf 235} (2001), no.~1, 243--260.

\bibitem{buch:grothendieck}
\bysame, \emph{Grothendieck classes of quiver varieties}, Duke Math. J.
{\bf 115} (2002), no. 1, 75--103.

\bibitem{buch:littlewood-richardson}
\bysame, \emph{A {L}ittlewood-{R}ichardson rule for the ${K}$-theory of
  {G}rassmannians}, Acta Math. {\bf 189} (2002), 37--78.

\bibitem{buch.fulton:chern}
A.~S. Buch and W.~Fulton, \emph{Chern class formulas for quiver
  varieties}, Invent.~Math. {\bf 135} (1999), 665-697.

\bibitem{schquiv}
A.~S. Buch, A.~Kresch, H.~Tamvakis and A.~Yong, \emph{Schubert
  polynomials
and quiver formulas}, Duke Math. J., to appear.

\bibitem{fomin.kirillov:grothendieck}
S.~Fomin and A.~N. Kirillov, \emph{Grothendieck polynomials and the 
{Y}ang-{B}axter equation}, Proceedings of the 6th Intern. Conf. on
 Formal Power Series and Algebraic Combinatorics, DIMACS (1994),
183--190.

\bibitem{fomin.kirillov:yang-baxter} 
\bysame, \emph{The {Y}ang-{B}axter equation, symmetric functions,
  and {S}chubert polynomials}, Proceedings of the 5th Conference on
  Formal Power Series and Algebraic Combinatorics (Florence, 1993),
  Discrete Math.  {\bf 153} (1996), no. 1-3, 123--143.

\bibitem{fulton:usp} W. Fulton,
{\em Universal Schubert polynomials},
Duke Math. J. {\bf 96} (1999), no. 3, 575--594.  

\bibitem{fulton.lascoux:pieri}
W.~Fulton and A.~Lascoux, \emph{A {P}ieri formula in the
  {G}rothendieck ring of a flag bundle}, Duke Math. J. \textbf{76}
  (1994), 711--729.

\bibitem{kirillov:cauchy} A. N. Kirillov, {\em Cauchy identities for
    universal Schubert polynomials}, Zap. Nauchn. Sem.  S.-Peterburg.
  Otdel. Mat. Inst. Steklov. (POMI) {\bf 283} (2001), 123--139.

\bibitem{lakshmibai.magyar:degeneracy}
V.~Lakshmibai and P.~Magyar, \emph{Degeneracy schemes, quiver schemes, and
  {S}chubert varieties}, Internat. Math. Res. Notices \textbf{1998}, no.~12,
  627--640. 
%\MR{99g:14065}

\bibitem{lascoux:grothendieck}
A.~Lascoux, \emph{Transition on Grothendieck polynomials},
Physics and combinatorics, 2000 (Nagoya),  164--179,
World Sci. Publishing, River Edge, NJ, 2001. 

\bibitem{lascoux.schutzenberger:schubpoly}
A.~Lascoux and M.-P.~Sch\"{u}tzenberger, \emph{Polyn\^{o}mes de
  Schubert}, C.R. Acad. Sci. Paris S\'{e}r. I Math. {\bf 294} (1982),
447--450.

\bibitem{lascoux.schutzenberger:grothendieck}
\bysame, \emph{Structure de Hopf de l'anneau de cohomologie et de
  l'anneau de Grothendieck d'une vari\'{e}t\'{e} de drapeaux},
C.R. Acad. Sci. Paris S\'{e}r. I Math. {\bf 295} (1982), 629--633.

\bibitem{Macdonald:schub_notes}
I.~G.~Macdonald, \emph{Notes on Schubert polynomials}, Publ. LACIM
{\bf 6}, Univ.~de Qu\'{e}bec \'{a} Montr\'{e}al, Montr\'{e}al, 1991.

\end{thebibliography}
%\input{bibliography.tex}

\end{document}